\documentclass[]{article} %% The same for {book}

\usepackage{graphicx, amssymb, amsmath}
\usepackage{hyperref}

\graphicspath{{./}{figures/}}

\newcommand{\rr}{\mathbb{R}}
\newcommand{\rn}{\mathbb{R}^n}
\newcommand{\rk}[1]{\mathbb{R}^{k_#1}}

\newtheorem{theo}{\bf Theorem}

\newtheorem{defi}[theo]{\bf Definition}

\begin{document}
\title{Numerical Approximation of Nash Equilibria for a Class of Non-Cooperative Differential Games}

\author
{Simone Cacace\thanks{
Dipartimento di Matematica, SAPIENZA - Universit\`a di Roma, p.le Aldo Moro, 5 - 00185 Rome, Italy. E-mail address: cacace@mat.uniroma1.it}, 
  Emiliano Cristiani\thanks{
  Dipartimento di Matematica, SAPIENZA - Universit\`a di Roma, p.le Aldo Moro, 5 - 00185 Rome, Italy. E-mail address: emiliano.cristiani@gmail.com},
  Maurizio Falcone\thanks{
  Dipartimento di Matematica, SAPIENZA - Universit\`a di Roma, p.le Aldo Moro, 5 - 00185 Rome, Italy. E-mail address: falcone@mat.uniroma1.it 
  (Corresponding author)}
}
\date{December 23, 2010}
\maketitle
\textbf{Abstract.} 
In this paper we propose a numerical method to obtain an approximation of Nash equilibria for multi-player non-cooperative games with a special structure. 
We consider the infinite horizon problem in a case which leads to a system of Hamilton-Jacobi equations. 
The numerical method is based on the Dynamic Programming Principle for every equation and on a global fixed point iteration. 
We present the numerical solutions of some two-player games in one and two dimensions.
The paper has an experimental nature, but some features and properties of the approximation scheme are discussed.


\section{Introduction}
The notion of Nash equilibrium is considered one of the most important achievements in the last century for its impact on the analysis of economic and social phenomena. 
It is worth to note that the formal definition of Nash equilibria  \cite{N51} has opened new research directions and has attracted the interest of several mathematicians to new fields of application. After the pioneering work of Von Neumann and Morgenstern \cite{VNM44}, the use of mathematically-based analysis in the study of economic sciences has received new impulse from the work of Nash.

There is extensive literature dealing with Nash equilibria for non-cooperative games, however the analysis of Nash equilibria for nonzero-sum multi-player differential games is more limited  (see e.g. the monographs \cite{BO89,F71} for an presentation of this theory). Moreover, only few papers  give a characterization of the value functions for the players in terms of partial differential equations, as it was the case for control problems and for zero-sum differential games described for example in \cite{BCD97}. More precisely, we know that under suitable regularity assumptions, if the value functions of a non-cooperative nonzero-sum multi-player differential game exist, they satisfy a system of first-order Hamilton-Jacobi equations, see \cite{AHbook}.
Typically theoretical results about these problems are valid only in very special cases and they are essentially limited to games in one dimension and with a simple dynamics, see e.g. \cite{B74,BP06,BS04,CP03}. More important, it is well known that the system of Hamilton-Jacobi equations can be ill-posed. To our knowledge, there are no theoretical results if the dimension of the problem is greater than one.

From the numerical point of view the situation  is even more deceiving since few results are available for Nash equilibria  in the framework of non-cooperative games. In particular, we mention the recent paper \cite{FFP09} where an approximation of (static) Nash equilibria based on Newton methods is presented, and the paper \cite{FK10} where the approximation is obtained via penalty methods.

Our goal is to construct an approximation scheme for Nash equilibria for non-cooperative differential games starting from the characterization obtained via Hamilton-Jacobi equations. Following \cite{BP06}, we deal with the system of stationary Hamilton-Jacobi equations giving a characterization of the value functions for a class of infinite-horizon games with nonlinear costs exponentially discounted in time. To this end we will extend to the system the class of dynamic programming approximation schemes studied for zero-sum differential games. The interested reader can find in \cite{F06} a detailed analysis of the schemes corresponding to pursuit-evasion games as well as some numerical tests (see also \cite{F97} for the origin of these methods and some control applications). At our best knowledge, the approximation scheme proposed in this paper is the first for Nash equilibria in the framework of differential games.

The paper is organized as follows.
In Section 2 we set up the problem, introduce the notations and recall the main results giving the characterization of the value functions. In Section 3 we introduce the semi-discrete and fully-discrete approximation schemes, and we describe the fixed point iterative scheme for the system of Hamilton-Jacobi equations. Some remarks about the fixed point algorithm are also discussed. Finally, in Section 4 we present the numerical results for some problems in dimension one and two.

\section{Setting the problem}
Let us consider a $m\hbox{-}$player non-cooperative differential game with controlled dynamics
\begin{equation}\label{dynamics}
\left\{\begin{array}{l}
 \dot y(t)=f\big(y(t),\alpha_1(t),...,\alpha_m(t)\big) \\
y(0)=x
\end{array}
\right. 
\end{equation}
where $t>0$, $x\in\rn$, $f:\rn\times\rk{1}\times ...\times\rk{m}\to\rn$ and $\displaystyle \alpha_i:[0,+\infty)\to A_i$  is the (open-loop) control associated to the $i$-th player ($i=1,...,m$) within a set of admissible control values $A_i\subseteq\rk{i}$, $k_i\geq 1$. We set 
$\mathcal{A}_i=\big\{\alpha_i:[0,+\infty)\to A_i\,,\,\alpha_i\,\mbox{ measurable}\big\}$ for $i=1,...,m$. In order to simplify the notations, we also set $\alpha(\cdot)=(\alpha_1(\cdot),\ldots,\alpha_m(\cdot))$ 
and we denote by $y_x(t;\alpha(\cdot))$ the corresponding solution of the Cauchy problem (\ref{dynamics}), i.e. the trajectory starting at $x$ which evolves according to  the strategies $\alpha(\cdot)$ of the $m$ players.\\
We consider the infinite horizon problem, where each player has a running cost discounted exponentially in time. More precisely, for $i=1,...,m$, we take $\lambda_i>0$, $\psi_i:\rn\times\rk{i}\times...\times\rk{m}\to\rr$ and 
we define the cost functionals
\begin{equation}\label{cost}
J_i(x,\alpha(\cdot))=\int_0^{+\infty}\hskip-15pt\psi_i\big(y_x(t;\alpha(\cdot)),\alpha(\cdot)\big)e^{-\lambda_i t}\,dt\,,\quad i=1,\dots, m. 
\end{equation}
We say that a $m$-tuple of feedback controls $a^*(y)=(a_1^*(y),...,a_m^*(y))$ (i.e. functions $a_i^*:\rn\to A_i$, $i=1,...,m$, depending on the state variable) is a \emph{Nash non-cooperative equilibrium solution} for the game (\ref{dynamics}) if
\begin{equation}\label{nash_eq}
J_i(x,a^*)=\hskip-5pt\min_{\alpha_i(\cdot)\in\mathcal{A}_i}\hskip-5pt J_i(x,a_1^*,...a_{i-1}^*,\alpha_i(\cdot),a_{i+1}^*,...,a_m^*)\,,\quad i=1,...,m\,. 
\end{equation}
Note that for every given feedback control $a_i:\rn\to A_i$, $i=1,...,m$, and every path $y(t)$, $t>0$, we can always define the corresponding open loop control as $a_i(y(t))
\in \mathcal{A}_i$, so that, in the above definition, $J_i(x,a^*)$ is the $i$-th cost associated to the trajectory $y_x(t;a^*)$ in which \emph{all} the optimal Nash strategies $a^*$ are implemented, namely the solution of
\begin{equation}\label{dynamics_opt}
\left\{\begin{array}{l}
 \dot y(t)=f\big(y(t),a^*(y(t))\big) \\
y(0)=x\,.
\end{array}
\right. 
\end{equation}
On the other hand the term $J_i(x,a_1^*,...a_{i-1}^*,\alpha_i(\cdot),a_{i+1}^*,...,a_m^*)$ in (\ref{nash_eq}) is the $i$-th cost associated to the trajectory $y_x(t;a_1^*,...,a_{i-1}^*,\alpha_i(\cdot),a_{i+1}^*,...,a_m^*)$ corresponding to the solution of
\begin{equation}\label{dynamics_opt-i}
\left\{\begin{array}{l}
 \dot y(t)=f\big(y(t),a_1^*(y(t)),...,a_{i-1}^*,\alpha_i(\cdot),a_{i+1}^*,...,a_m(y(t))\big) \\
y(0)=x
\end{array}
\right. 
\end{equation}
where \emph{only} the strategy $\alpha_i(\cdot)$ is chosen in $\mathcal{A}_i$.\\
The definition of Nash equilibrium means that if the $i$th player replaces his optimal control $a_i^*$ with any other strategy $\alpha_i(\cdot)\in\mathcal{A}_i$, then his running cost $J_i$ increases, 
assuming that the remaining players keep their own controls frozen. The $m\hbox{-}$tuple $a^*$ is then optimal in the sense that no player can do better for himself, since he can not cooperate with any other player.

Let us assume that such a Nash equilibrium $a^*$ exists for our game problem. 
%\note{aggiungere}{\bf see \cite{??} for sufficient conditions and some examples?? trovare una referenza precisa, Friedman ??}\\ 
For $i=1,...,m$ we define the value function $u_i:\rn\to\rr$ as the minimal cost $J_i$ associated to $a^*$. More precisely, for every $x\in\rn$ and $i=1,...,m$, we set
\begin{equation}\label{value_function}
 u_i(x)=J_i(x,a^*)\,.
\end{equation}
 Then, it can be proved that all the $u_i$ satisfy a Dynamic Programming Principle and by standard arguments we can derive a system of Hamilton-Jacobi equations for $u_1,...,u_m$, in which the feedback control $a^*(x)=(a_1^*(x),...,a_m^*(x))$ depends on the state variable $x$ also through the gradients $\nabla u_1(x),...,\nabla u_m(x)$. We get
\begin{equation}\label{HJsystem}
\lambda_i u_i(x)=H_i\big(x,\nabla u_1(x),...,\nabla u_m(x)\big)\qquad x\in\rn\,,i=1,...,m\,,
\end{equation}
where, for every  $i=1,...,m$ and for every $x$, $p_1,...,p_m\in\rn$, the Hamiltonians $H_i:\rr^{n+nm}\to\rr$ are given by
\begin{equation}\label{hamiltonians}
H_i(x,p_1,...,p_m)=p_i\cdot f\big(x,a_1^*(x,p_1,...,p_m),...,a_m^*(x,p_1,...,p_m)\big)+
\end{equation}
$$+\psi_i\big(x,a_1^*(x,p_1,...,p_m),...,a_m^*(x,p_1,...,p_m)\big)\,.
$$
Moreover, for every $x\in\rn$ and $i=1,...,m$, the following property holds:
\begin{eqnarray}\label{hamiltonians_min}
 H_i\big(x,\nabla u_1(x),...,\nabla u_m(x)\big)=\hspace{5cm}\nonumber\\
=\min_{a_i\in A_i}\Big\{\nabla u_i(x)\cdot f\big(x,a_1^*(x),...,a_{i-1}^*(x),a_i,a_{i+1}^*(x),...,a_m^*(x)\big)\,+ 
\\ \nonumber
  +\,\psi_i\big(x,a_1^*(x),...,a_{i-1}^*(x),a_i,a_{i+1}^*(x),...,a_m^*(x)\big)\Big\}\,.\hspace{2.3cm}\nonumber\
\end{eqnarray}
We remark that in (\ref{hamiltonians_min}) the minimum is performed among all the control values $a_i\in A_i$ and not among all the strategies $\alpha_i(\cdot)\in\mathcal{A}_i$.

\section{Numerical approximation}
This section is devoted to the introduction of the semi-discrete and the fully-discrete scheme for the system of Hamilton-Jacobi equations (\ref{HJsystem}), the fixed-point algorithm, and to a brief discussion about the algorithm.
\subsection{Semi-discrete and fully-discrete scheme}\label{fully_discrete}
Here we propose a numerical scheme to compute the value functions $u_i$, $i=1,\dots, m$ defined in (\ref{value_function}). 
To simplify the presentation and for computational purposes, we deal with a two-player game ($m=2$), with scalar controls (i.e. $k_1=k_2=1$) and a dynamics in one or two dimensions ($n=1,2$).
 Moreover,  we assume the discount rates $\lambda_1=\lambda_2=1$.

In order to discretize the Hamiltonians $H_i$ we use a semi-Lagrangian scheme, which usually gives more accurate results than finite difference schemes or others schemes. The reader can find in \cite{F97} 
a comprehensive introduction to this subject. As usual we first obtain a semi-discrete scheme, introducing a fictitious time step $h>0$, which reads as
\begin{equation}\label{SLsemidiscrete}
\left\{
\begin{array}{l}
u_1(x)=\displaystyle\min_{a_1\in A_1}\left\{ \frac{1}{1+h} u_1(x+h\,f(x,a_1,a_2^*)) + \frac{h}{1+h}\psi_1(x,a_1,a_2^*) \right\} \\ \\
u_2(x)=\displaystyle\min_{a_2\in A_2}\left\{ \frac{1}{1+h} u_2(x+h\,f(x,a_1^*,a_2)) + \frac{h}{1+h}\psi_2(x,a_1^*,a_2) \right\}
\end{array}
\right.
\end{equation}
where we remind that the control $a^*=(a_1^*,a_2^*)$ depends on $x$ and $\nabla u_1(x)$, $\nabla u_2(x)$. Note that the time step $h$ can be interpreted as the discrete time corresponding to the approximation
of the trajectories starting from $x$ and moving according to the dynamics $f$.

Now let us consider a subdomain $\Omega\subset\rn$ in which we look for the approximate solutions of the system described above. We discretize $\Omega$ by means of a uniform grid of size $\Delta x$ denoted by $G=\{x_1,\ldots,x_N\}$, where $N$ is total number of nodes. For $i=1,2$ we denote by $U_i\in\rr^N$ the vector containing the values of $u_i$ at the grid nodes, \emph{i.e.} $(U_i)_j=u_i(x_j)$, $j=1,\ldots,N$. Moreover, for $i=1,2$, let $A_i^\#$ be a finite discretization of the set of admissible controls $A_i$. 

Note that  even for $x=x_j$ the point $z(x_j,a^*)=x_j+h\,f(x_j,a^*))$ appearing in (\ref{SLsemidiscrete}) will not coincide with a node of the grid $G$. Then, an interpolation is needed to compute the value of $u_i$ at $z$ (this is the main difference with respect to a standard finite difference approximation).
 First of all, denoting by $\displaystyle\|f\|_\infty$ the infinity norm of $f$ with respect to all its variables and choosing  
\begin{equation}\label{h}
h=\frac{\Delta x}{\|f\|_\infty}
\end{equation} 
we guarantee that the point $z$ lies in one of the first neighbouring cells around $x_j$, denoted by $I(z)$. Then, the value $u_i(z)$ is reconstructed by a linear interpolation of the values of $u_i$ at the vertices of the cell $I(z)$ (see \cite{F97} for more details and \cite{CFF04} for an efficient algorithm in high dimension). Let $\Lambda(a)$ denote the $N\times N$ matrix which defines the interpolation: for $j=1,\ldots,N$, the $j$th row contains the weights to be given to the values $(U_i)_1,\ldots,(U_i)_N$ in order to compute $u_i(z(x_j,a))$. The fully-discrete version of the system (\ref{SLsemidiscrete}) can be written in compact form as
\begin{equation}\label{SLfullydiscrete}
\left\{
\begin{array}{l}
(U_1)_j=\displaystyle\min_{a_1\in A_1^\#}\left\{ \frac{1}{1+h} \left( \Lambda(a_1,a_2^*)\ U_1\right)_j + \frac{h}{1+h}\psi_1(x_j,a_1,a_2^*) \right\} \\ \\
(U_2)_j=\displaystyle\min_{a_2\in A_2^\#}\left\{ \frac{1}{1+h} \left( \Lambda(a_1^*,a_2)\ U_2\right)_j + \frac{h}{1+h}\psi_2(x_j,a_1^*,a_2) \right\}.
\end{array}
\right.
\end{equation}
where the index $j$ varies from 1 to $N$.\\ 
We remark that for an actual implementation, the values of $U$ at the boundary nodes of the grid must be managed apart, in order to impose boundary conditions on $\partial \Omega$. 
In Section \ref{tests} 
dedicated to numerical tests, we will use Dirichlet boundary conditions, taking notice of the influence they could have on the numerical solution inside $\Omega$.\\  
We are ready to describe the algorithm we actually implemented.
\vskip0.3cm
\emph{Fixed point algorithm}
\begin{enumerate}
\item Choose two tolerances $\varepsilon_1>0$ and $\varepsilon_2>0$. \\
Set $k=0$. Choose an initial guess for the values of $U_i$, $i=1,2$ and denote them by $U_i^{(0)}$.
\item For $j=1,\ldots,N$

\begin{enumerate} 
%\note{bisogna dire qualcosa su possibili alternative}
\item Find $a^*\in A_1^\# \times A_2^\#$ such that
$$%\begin{equation}\label{SLfullydiscrete}
\left\{
\begin{array}{l}
a_1^*=\displaystyle\arg\!\min\limits_{\!\!\!\!\!\!a_1\in A_1^\#}\left\{ \frac{1}{1+h} \left( \Lambda(a_1,a_2^*)\ U_1^{(k)}\right)_j + \frac{h}{1+h}\psi_1(x_j,a_1,a_2^*) \right\} \\ \\
a_2^*=\displaystyle\arg\!\min\limits_{\!\!\!\!\!\!a_2\in A_2^\#}\left\{ \frac{1}{1+h} \left( \Lambda(a_1^*,a_2)\ U_2^{(k)}\right)_j + \frac{h}{1+h}\psi_2(x_j,a_1^*,a_2) \right\}.
\end{array}
\right.
$$%\end{equation}
Note that the search for $a^*$ can be done in an exhaustive way, due to the  fact that the set $A_1^\# \times A_2^\#$ is finite.
\item If $a^*$ is found, go to Step $(c)$, otherwise stop (if more than one Nash optimal control is available, we select the first one we find).
\item Solve, for $i=1,2$,
\begin{equation}\label{SLscheme_algorithm}
(U_i^{(k+1)})_j=\frac{1}{1+h} \left( \Lambda(a^*)\ U_i^{(k)}\right)_j + \frac{h}{1+h}\psi_i(x_j,a^*)
\end{equation}
\end{enumerate}
\item If $\|U_i^{(k+1)}-U_i^{(k)}\|_{\rr^N}<\varepsilon_i$, $i=1,2$, then stop, else go to step 2 with $k \leftarrow k+1$.
\end{enumerate}
\subsection{Some remarks about the fixed point algorithm}
We introduce a vector $U$ containing the two value functions,
$$U=((U_1)_1,\ldots,(U_1)_N,(U_2)_1,\ldots,(U_2)_N)\in\rr^{2N}$$
and a vector $\Psi$ containing the two cost functions, 
$$\Psi(a)=(\psi_1(x_1,a),\ldots,\psi_1(x_N,a),\psi_2(x_1,a),\ldots,\psi_2(x_N,a))\in\rr^{2N}\,.$$
We also define a fixed point operator $F=(F_1,\ldots,F_{2N}):\rr^{2N}\rightarrow \rr^{2N}$, given component-wise by
\begin{equation}\label{defF}
F_j(U)=\frac{1}{1+h} \left( \Lambda(a^*)\ U\right)_j + \frac{h}{1+h}(\Psi(a^*))_j\,,\qquad j=1,\ldots,2N\,,
\end{equation} 
where, with an abuse of notation, we denoted again by $\Lambda(a)$ the block matrix
$$
\left(
\begin{array}{c|c}
\Lambda(a) & 0 \\
\hline
0 &  \Lambda(a)
\end{array}
\right).
$$
In this way the scheme (\ref{SLscheme_algorithm}) can be written as $U=F(U)$.

In the case of a one-player game or two-player zero-sum games with minmax equilibrium, it can be easily proven that the corresponding fixed point operator (defined similarly as before) is a contraction map, 
see for example \cite{F06}. More precisely, denoted by $G$ the operator, it can be proved that, for some norm $\|\cdot\|$,
\begin{equation}\label{Gcontrazione}
\|G(U)-G(V)\|\leq \frac{1}{1+h}\|U-V\|\,,\quad \forall ~U,V.
\end{equation}
In the case of Nash equilibria, the arguments leading to the proof of (\ref{Gcontrazione}) are not valid. 
Nevertheless, we could in principle compute the Jacobian matrix $J_F(U)$ of $F(U)$, i.e.,
\begin{equation}\label{JFU}
J_F(U)=\left(\frac{\partial F_q(U)}{\partial U_r}\right)_{q,r=1,\ldots,2N}
\end{equation}
and then study the behaviour of its norm. In this respect, it is important to note that the optimal control $a^*$ which appears in (\ref{defF}) actually depends on $U$, and then $\Lambda(a^*)$ depends on $U$. 
This makes the analytical computation of $J_F(U)$ extremely difficult, since we should know the derivative of $a^*$ with respect to $U$. Nevertheless, in some particular cases the computation can be easy. 
Indeed, since we assume that $A_1^\#\times A_2^\#$ is finite, it is reasonable to expect that small variation of $U$ do not produce a change of $a^*$, \emph{i.e.} the function $a^*(U)$ is ''piecewise constant''. 
It is worth to note that in the zones of $\rr^{2N}$ where $a^*(U)$ is constant,  $\|J_F(U)\|_\infty$ is easy computed and we have
\begin{equation}\label{normJFU}
\|J_F(U)\|_\infty=\frac{1}{1+h}<1
\end{equation}
as in the one-player case or the two-player zero-sum minmax case \cite{F06} (the result easily follows from the fact that the sum of the elements of each row of $\Lambda$ is 1, 
by definition of linear interpolation). Summarizing, it results that $F(U)$ can be piecewise contractive in the $\infty$-norm, according to the following
\begin{defi}\label{def:pc}
A function $ F: \rr^{2N} \rightarrow \rr^{2N}$ is {\em piecewise contractive} in the norm $\| \cdot \|$ if and only if :\\
a) $F$ is differentiable a.e. in $\rr^{2N}$;\\
b) $\|J_F(x)\| <1$ for all the points $x\in\rr^{2N}$ where $F$ is differentiable.
\end{defi}
From the numerical point of view, we can compute the Jacobian matrix $J_F(U)$ of $F(U)$
and check if $\|J_F(U)\|_\infty<1$ for some $U$. Such a computation can give some information about the behaviour of $F$ near a tentative solution. Moreover, we can plot the graph of some component $F_q$ of $F$ as a function of one of its variables $U_r$, again in order to understand the regularity of $F$.

\section{Numerical tests}\label{tests}
The aim of this section is to provide some numerical simulations in dimensions $n=1,2$ in order to give some insights about the well-posedness of system (\ref{HJsystem}). 
We believe that such results can be useful to lead future theoretical investigations.\\
Let us start by the one dimensional examples analyzed in \cite{BP06}. 
The game takes place in $\rr$ with scalar controls $a_1,\,a_2\in\rr$ and it is driven 
by the dynamics 
\begin{equation}\label{dynbressan}
f(x,a_1,a_2)=a_1+a_2\,.
\end{equation}
The cost functions are 
\begin{equation}\label{costbressan}
\psi_i(x,a_1,a_2)=h_i(x)+\frac{a_i^2}{2}\,,\qquad i=1,2\,,
\end{equation}
where the $h_i$'s are smooth functions satisfying $|h_i^\prime(x)|\leq C$ for some constant $C$. 

For reader's convenience, we recall here a notion of {\em admissible solution} for the system (\ref{HJsystem}) introduced in \cite{BP06}, 
for which an existence and uniqueness result is available (see \cite{BP06}, Theorem 4).  
\begin{defi}\label{admissible}
A function $u=(u_1,u_2):\rr\to\rr^2$ is called an {\em admissible solution} to the system (\ref{HJsystem}) if the following holds:
\begin{itemize}
 \item[(A1)] $u$ is absolutely continuous and satisfies system (\ref{HJsystem}) at a.e. point $x\in\rr$\,.
\item[(A2)] $u$ has sublinear growth at infinity, i.e., there exists a positive constant $C$ such that, for all $x\in\rr$,
$$|u(x)|\leq C(1+|x|)\,.$$ 
\item[(A3)] At every point $x\in\rr$, the derivative $u^\prime$ admits left and right limits $u^\prime(x-)$, $u^\prime(x+)$. At points where $u^\prime$ is discontinuous, 
these limits satisfy at least one of the conditions
$$u^\prime_1(x-)+u^\prime_2(x-)\geq 0\qquad\mbox{or}\qquad u^\prime_1(x+)+u^\prime_2(x+)\leq 0\,.$$
\end{itemize}
\end{defi}
 
\noindent{\bf Test 1}\\
We choose the dynamics (\ref{dynbressan}) and cost functions (\ref{costbressan}) with $h_i(x)\equiv0$, $i=1,2$. 
The corresponding system of Hamilton-Jacobi equations is 
\begin{equation}\label{systembressan}
\left\{
\begin{array}{l}
u_1=-\displaystyle\left(\frac{ u_1^\prime}{2} +u_2^\prime\right)u_1^\prime\,,\\ 
u_2=-\displaystyle\left(\frac{ u_2^\prime}{2} +u_1^\prime\right) u_2^\prime
\end{array}
\right.
\end{equation}
and admits at least the following three formal solutions:
$$
u_1(x)\equiv u_2(x)\equiv 0\,,
$$
$$
\overline{u}_1(x)=\left\{\begin{array}{ll} 0 & \mbox{if } |x|\geq 1 \\ -\frac{1}{2}(1-|x|)^2 & \mbox{otherwise\,,}\end{array}\right. \qquad \overline{u}_2(x)=0\,,
$$
and
$$
\widehat{u}_1(x)=-\frac{1}{2}x^2\,,\qquad \widehat{u}_2(x)=0\,,
$$
with optimal controls $a^*=(-u_1^\prime,-u_2^\prime)$, $\overline{a}^*=(-\overline{u}_1^\prime,-\overline{u}_2^\prime)$ and $\widehat a^*=(-\widehat u_1^\prime,-\widehat u_2^\prime)$ respectively.\\
%\note{calcolare e inserire!!!}
It is easy to check that $u=(u_1,u_2)$ is the only admissible solution among the three above. Indeed, the solution $\overline{u}$ does not satisfy condition (A3) in Definition \ref{admissible}, whereas the 
solution $\widehat{u}$ does not satisfy condition (A2) in Definition \ref{admissible}. %Moreover, the optimal control in feedback form for the admissible solution $u$ is given by $a^*=(-u_1^\prime,-u_2^\prime)$.

\begin{figure}[h!]
\vskip-10pt
\begin{tabular}{lr}
\hskip-25pt\includegraphics[width=.56\textwidth]{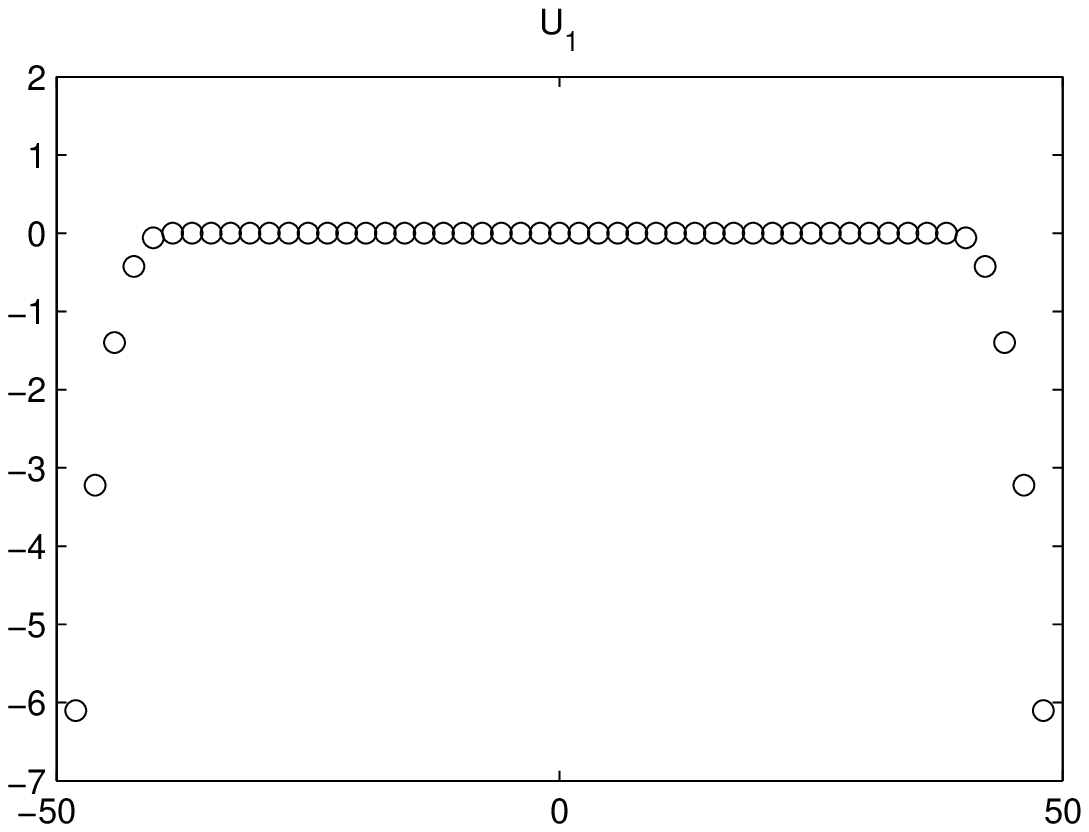} 
\hskip-20pt\includegraphics[width=.56\textwidth]{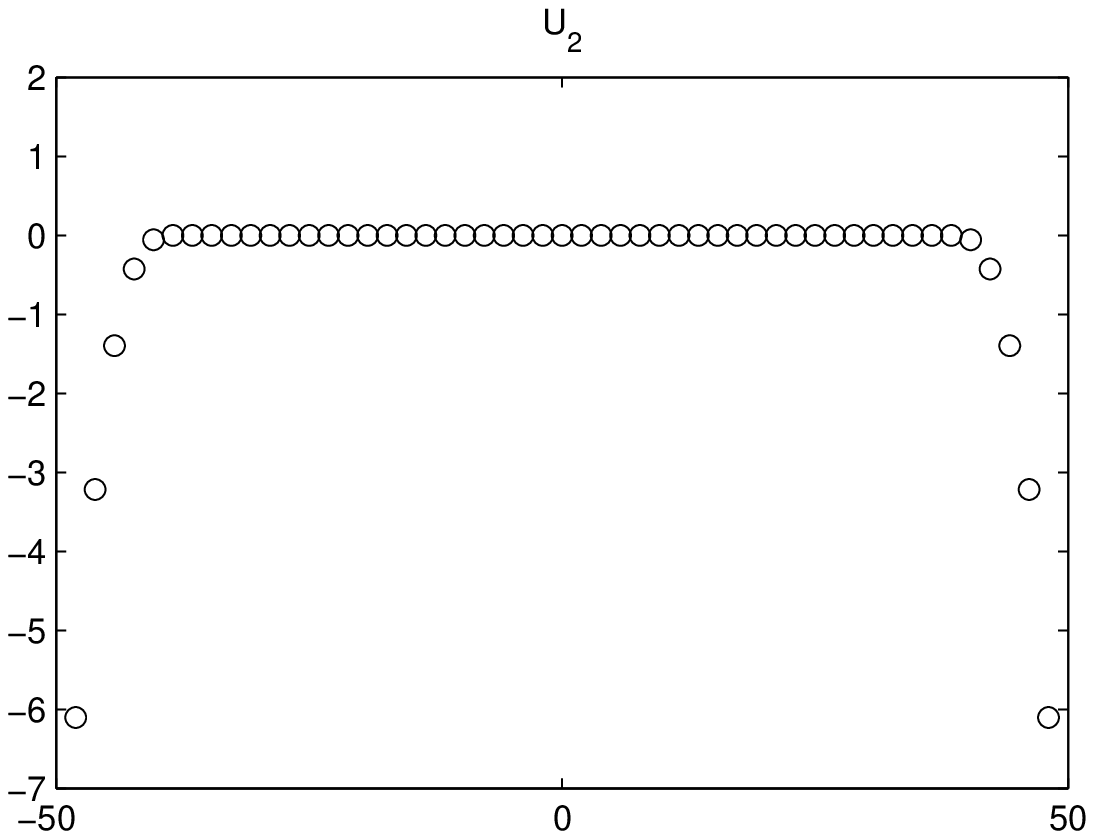}
\end{tabular}
\caption{Test 1. A boundary layer effect for the value functions $U_1$ (left) and $U_2$ (right).}
\label{T1_boundlayer}
\end{figure}

We solve system (\ref{systembressan}) by means of the semi-Lagrangian scheme presented in Section \ref{fully_discrete}. Let the numerical domain $\Omega=[-50,50]$ be discretized in $51$ nodes.
Moreover we approximate the sets of admissible controls $A_1=A_2=\rr$ with $101$ values uniformly distributed in $[-10,10]$. For $i=1,2$ we set the initial guess $U_i^{(0)}$ equal to a large constant. 
In this case a Nash equilibrium is always found and the scheme converges to the admissible solution $u\equiv 0$. We observe that the scheme converges to this solution also starting from the solution $u$ itself or from small perturbations of it.
Conversely, the two other solutions seem to be unstable: even if we choose them as initial guess, the scheme still converges to $u$. We only remark that the case of $\widehat{u}$ requires a larger set of controls, 
due to the large values attained in $\Omega$ by the optimal control $\widehat{a}^*=(-\widehat{u}_1^\prime,-\widehat{u}_2^\prime)=(x,0)$.\\ 
Concerning the Dirichlet boundary condition, we have noted that a large constant does not affect the solution inside the numerical domain and the scheme converges quite fast, 
whereas a negative value leads to the formation of boundary layers that dramatically slow down (and sometimes hinder) the convergence. For example, 
if we choose the value $-10$ as Dirichlet boundary condition for both $U_1$ and $U_2$, we obtain, after $200$ iterations, the results shown in Fig.\ \ref{T1_boundlayer}. \\

\begin{figure}[h!]
\vskip-10pt
\begin{tabular}{lr}
\hskip-25pt\includegraphics[width=.56\textwidth]{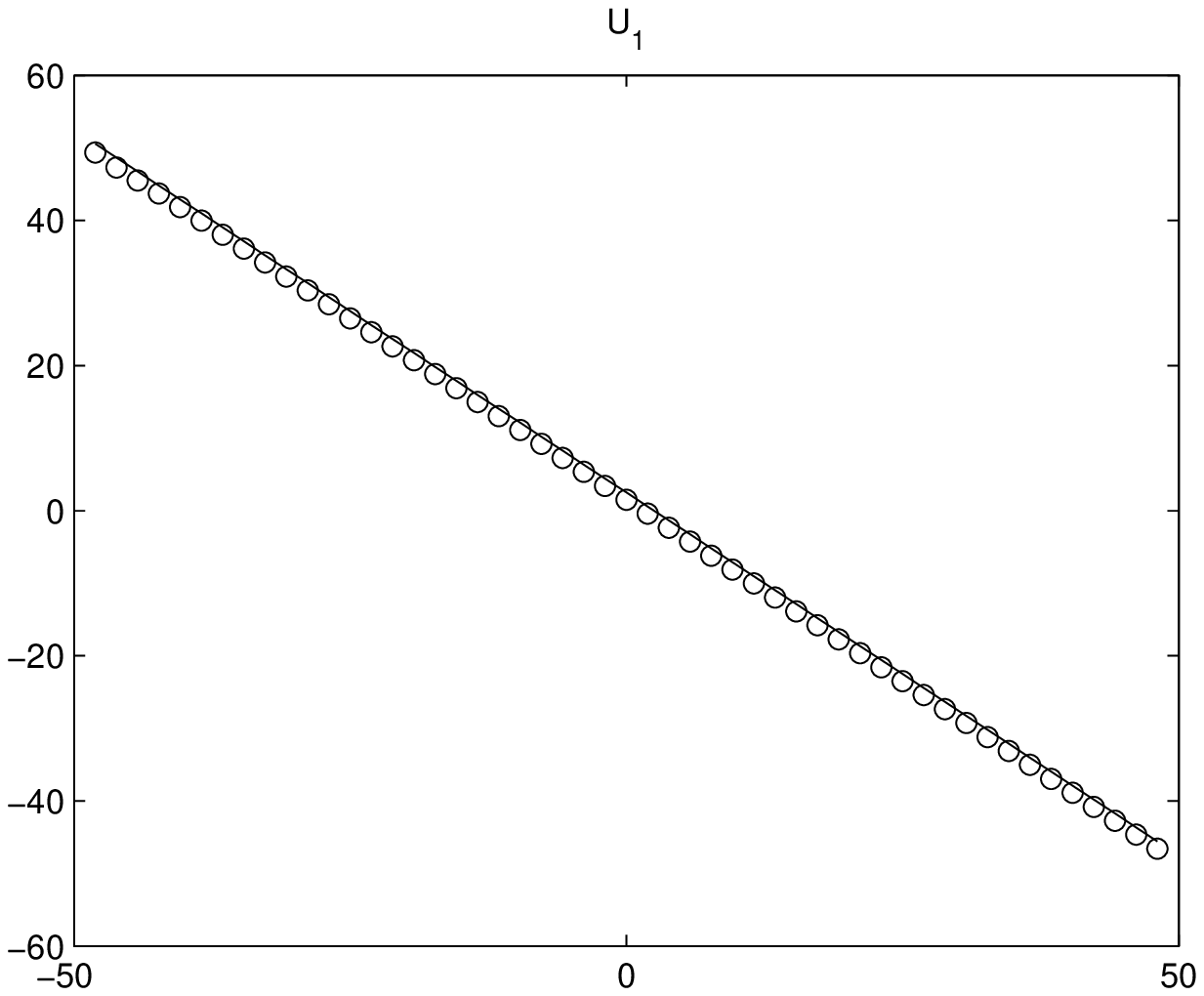} 
\hskip-20pt\includegraphics[width=.56\textwidth]{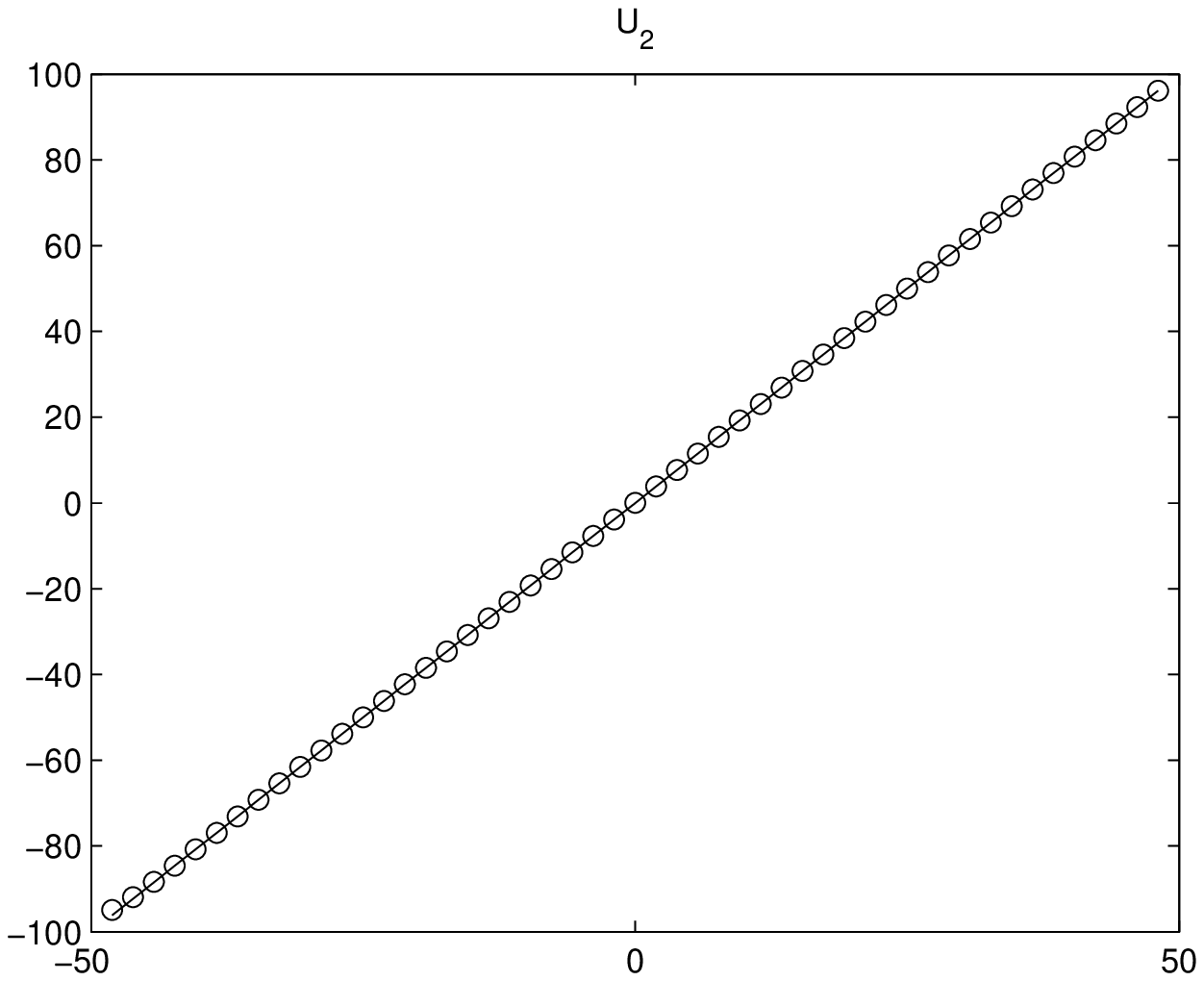}
\end{tabular}
\caption{Test 2. Comparison between exact solution (solid line) and approximate solution (circles). $u_1$ and $U_1$ are plotted on the left, $u_2$ and $U_2$ are plotted on the right.}
\label{T2}
\end{figure}

\begin{figure}[h!]
\begin{tabular}{lr}
\hskip-25pt\includegraphics[width=.56\textwidth]{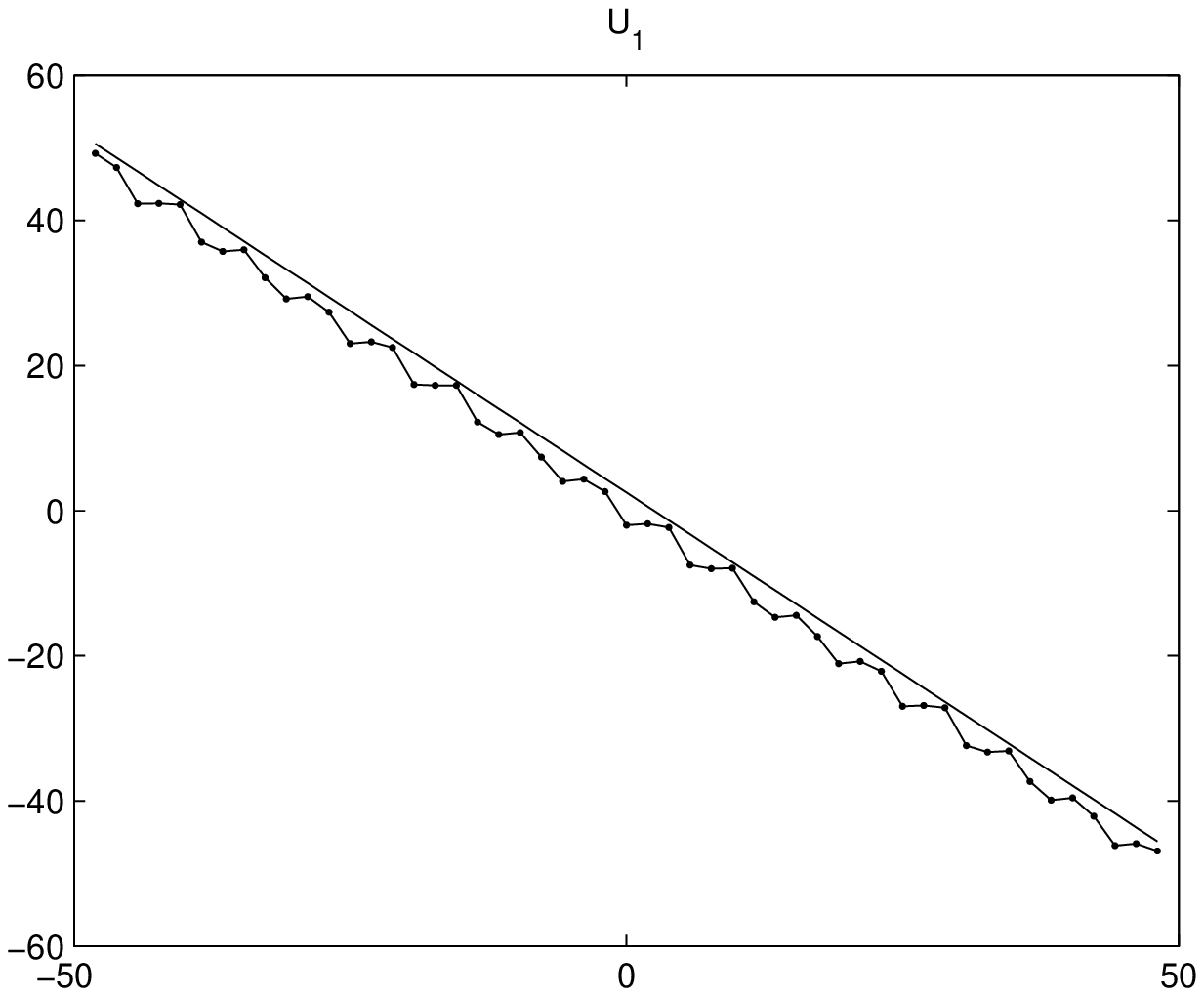} 
\hskip-20pt\includegraphics[width=.56\textwidth]{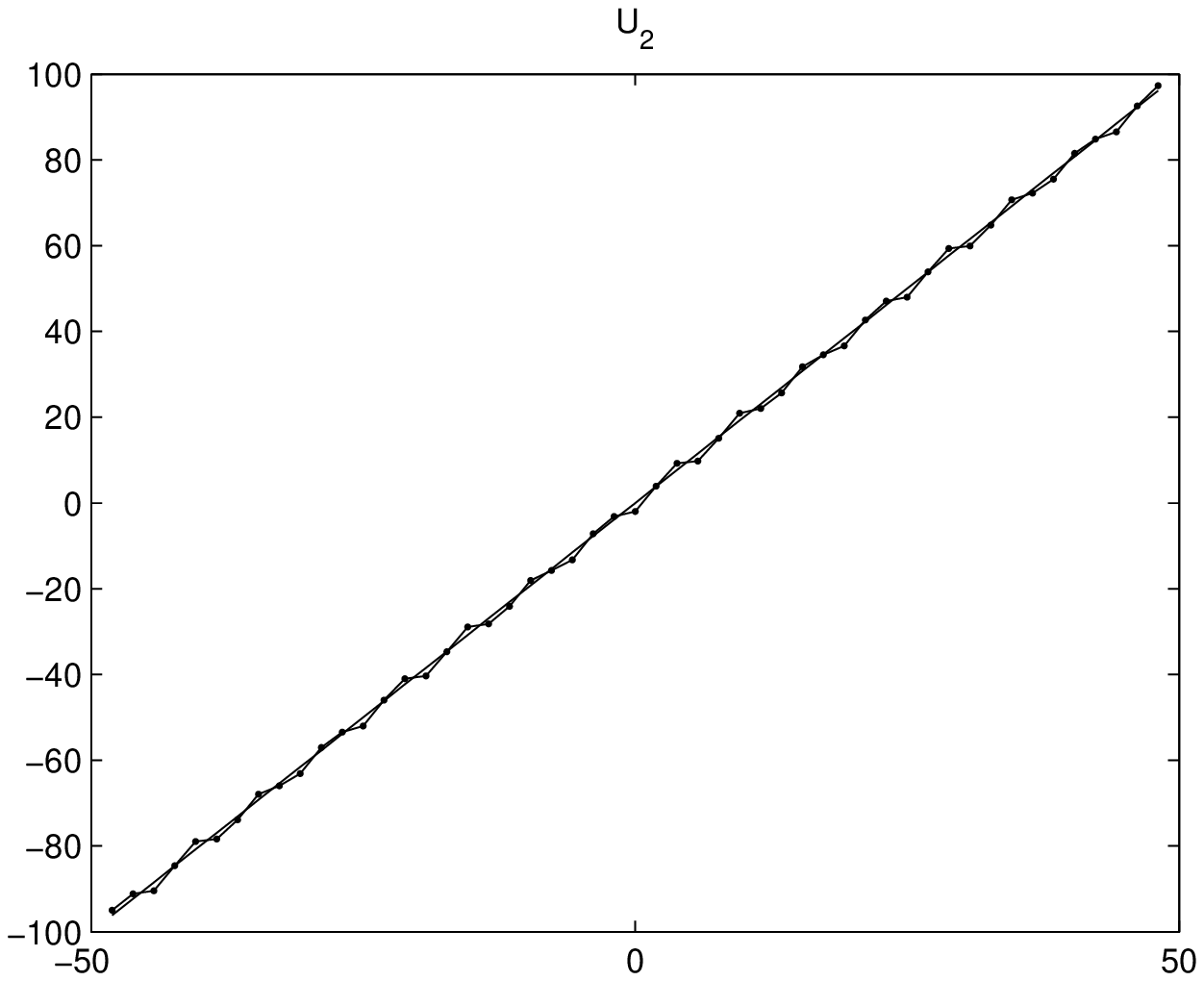}
\end{tabular}
\caption{Test 2. Comparison between exact solution obtained with $\delta=0$ (solid line) and approximate solution with $\delta=2$ (dots and line). $u_1$ and $U_1$ are plotted on the left, $u_2$ and $U_2$ are plotted on the right.}
\label{T2_ondine}
\end{figure}

\noindent{\bf Test 2}\\
Here we consider again the two-player game described by the dynamics (\ref{dynbressan}) and cost functions (\ref{costbressan}), where now we choose $h_i(x)=k_i x$ for $i=1,2$, with constants 
$k_1,k_2$ satisfying $k_1<0<k_2$ and $k_1+k_2\neq 0$. In this special case the solution is explicit, given by
$$
u(x)=\big(u_1(x), u_2(x)\big)=\big( k_1 x-k_1 k_2 -k_1^2/2\,,\,k_2 x-k_1 k_2 -k_2^2/2\big)\,.
$$
Moreover it is the unique admissible solution of system (\ref{HJsystem}), by virtue of Theorem $4$ proved in \cite{BP06}. 
We choose the same parameters of Test 1 and $k_1=-1$, $k_2=2$. 
Starting with the projection of $u$ on the grid, the algorithm reaches convergence immediately. Moreover, we observe convergence also starting from small perturbations of $u$ and from a constant initial guess $U^{(0)}\equiv C$ with 
$$C\geq \max_{x\in [-50,50]} \max\{u_1(x),u_2(x)\}\,.$$  
However convergence is not guaranteed for any initial guess.\\ 
Fig.\ \ref{T2} shows the approximate value functions compared with the exact solution corresponding to the choice $U^{(0)}\equiv C=150$.
\noindent Following \cite{BP06}, the existence and uniqueness result for admissible solutions holds in the more general case in which $h_i(x)$ is a sufficiently small perturbation of $k_i x$, for $i=1,2$. 
More precisely, for $k_1,k_2$ such that $k_1<0<k_2$ and $k_1+k_2\neq 0$, it is assumed that 
$$
|h_1^\prime(x)-k_1|\leq \delta\,,\qquad|h_2^\prime(x)-k_2|\leq \delta\,,
$$
for all $x\in\rr$ and a sufficiently small $\delta>0$.\\ 
We apply our algorithm also in this case, by choosing $\delta=2$ and
$$
h_1(x)=k_1 x -\delta \cos(x)\,, \qquad h_2(x)=k_2 x -\delta \cos(x)\,.
$$     
We choose a constant initial guess equal to 100 and the set of admissible controls is $[-300,300]$. The algorithm converges to a reasonable solution shown in Fig.\ \ref{T2_ondine}.\\

\begin{figure}[h!]
\begin{tabular}{lr}
\hskip-25pt\includegraphics[width=.56\textwidth]{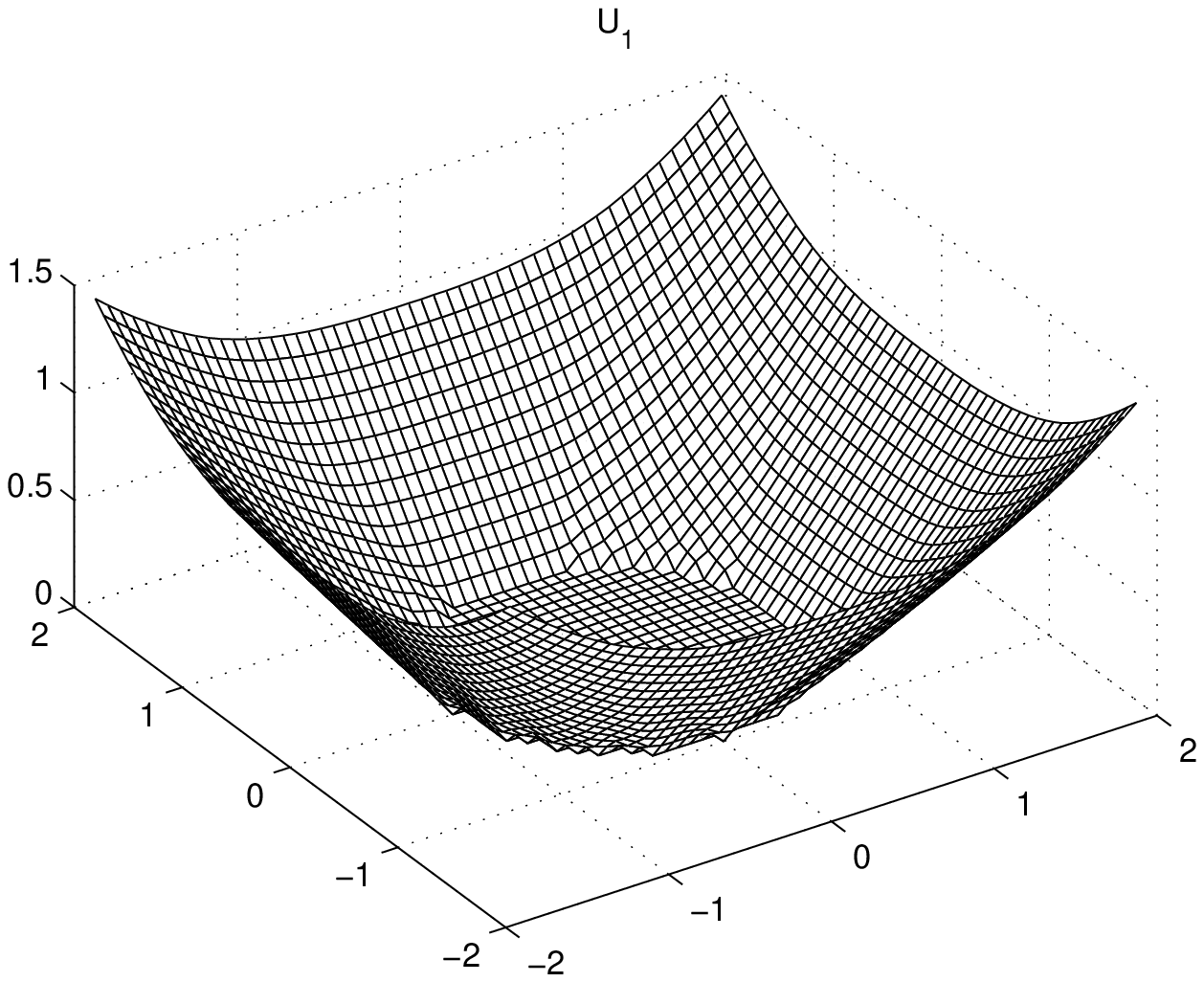} 
\hskip-20pt\includegraphics[width=.56\textwidth]{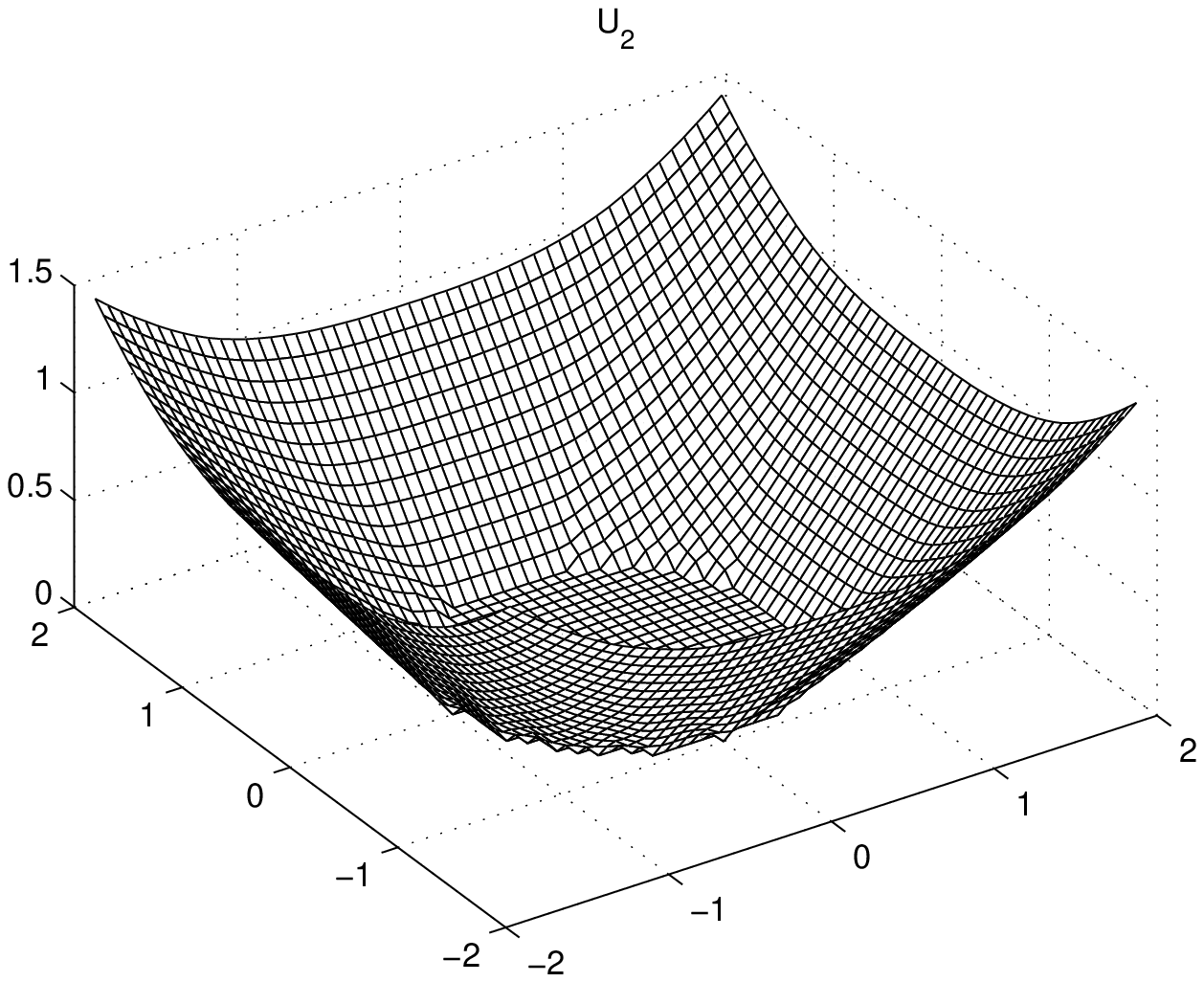}
\end{tabular}
\caption{Test 3. Approximate solutions $U_1$ (left) and $U_2$ (right).}
\label{T3}
\end{figure}

\begin{figure}[h!]
\begin{tabular}{lr}
\hskip-25pt\includegraphics[width=.56\textwidth]{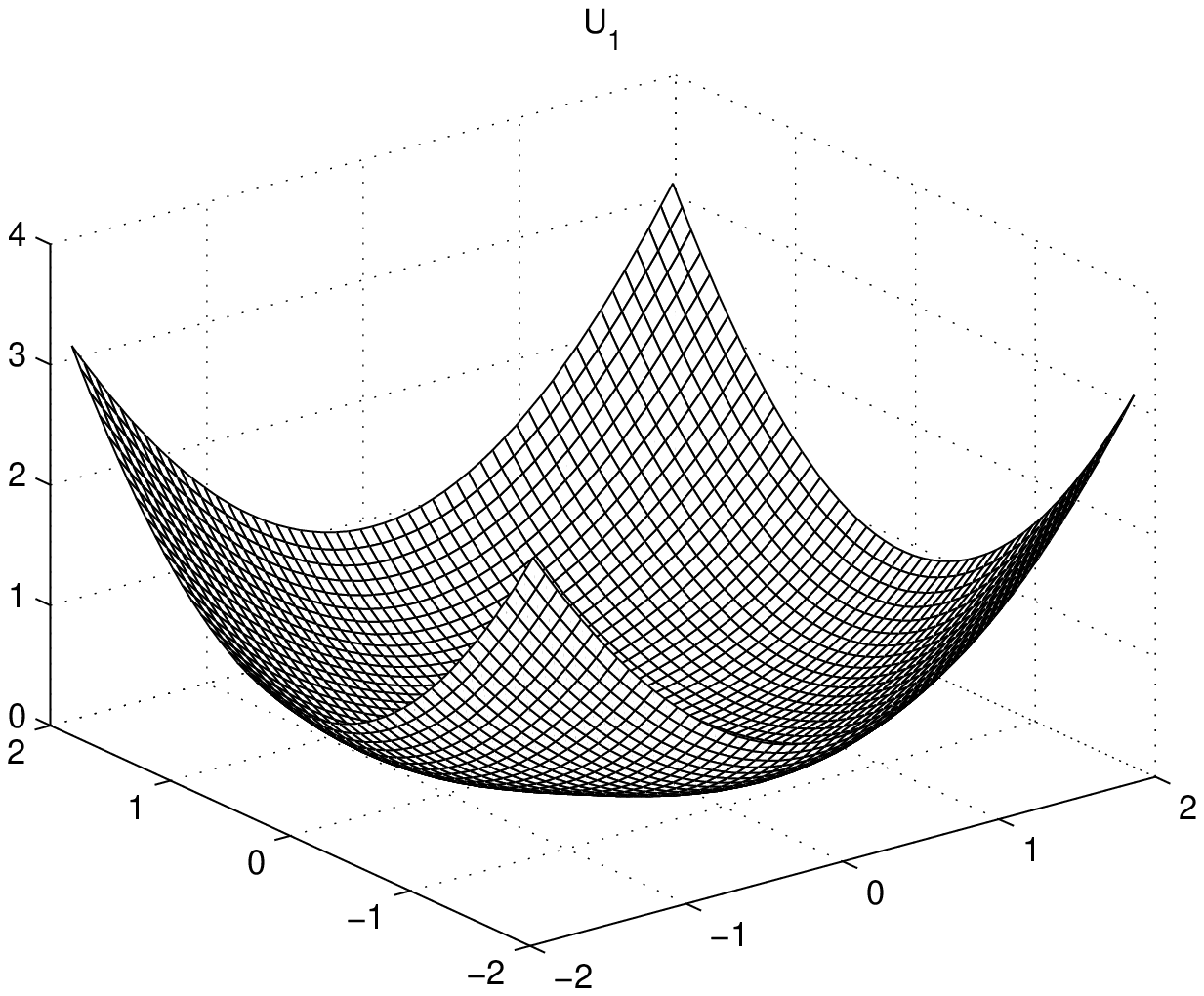} 
\hskip-20pt\includegraphics[width=.56\textwidth]{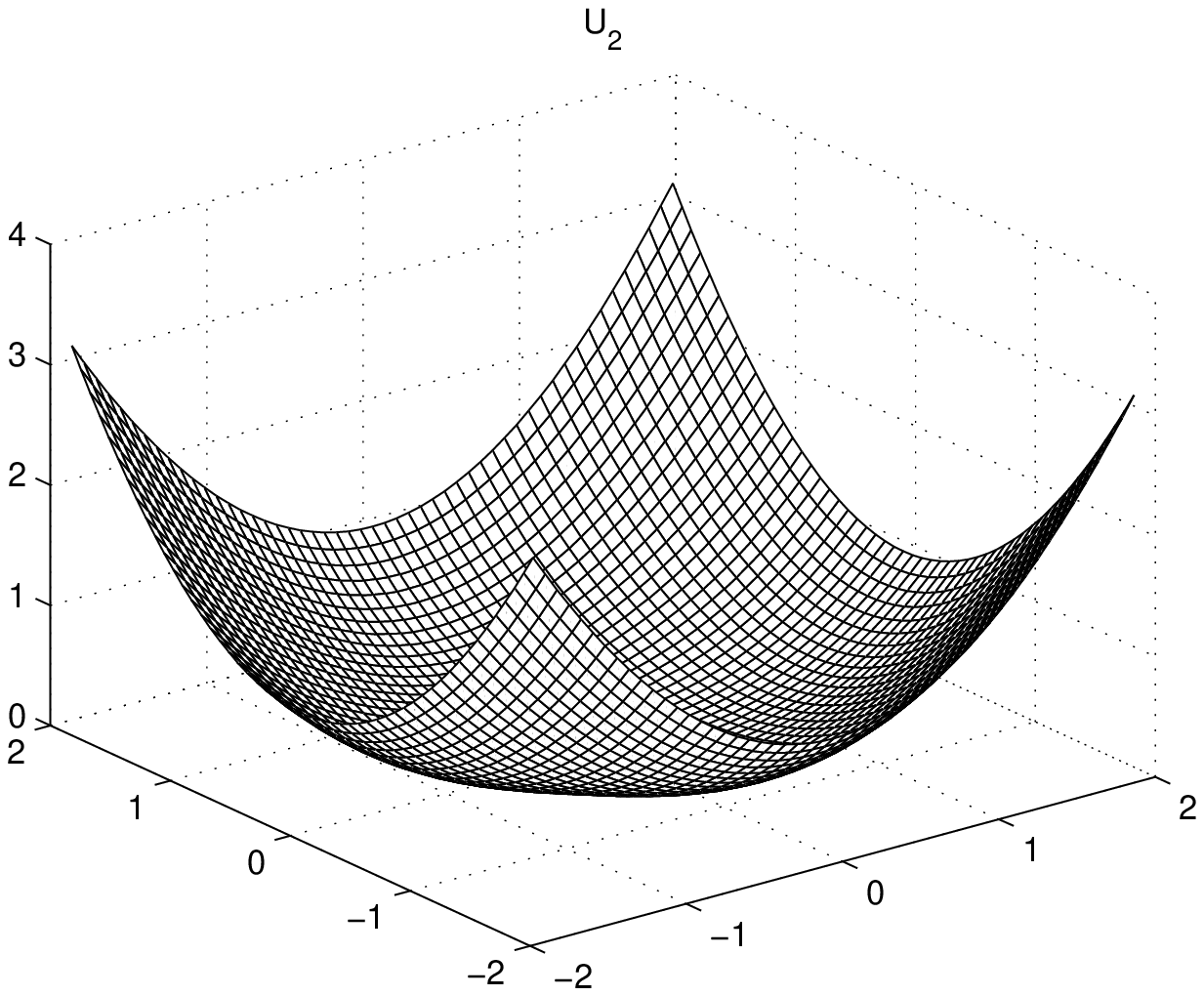}
\end{tabular}
\caption{Test 4. Approximate solutions $U_1$ (left) and $U_2$ (right) after 1000 iterations.}
\label{T4_U}
\end{figure}

\begin{figure}[h!]
\centering\includegraphics[width=.8\textwidth]{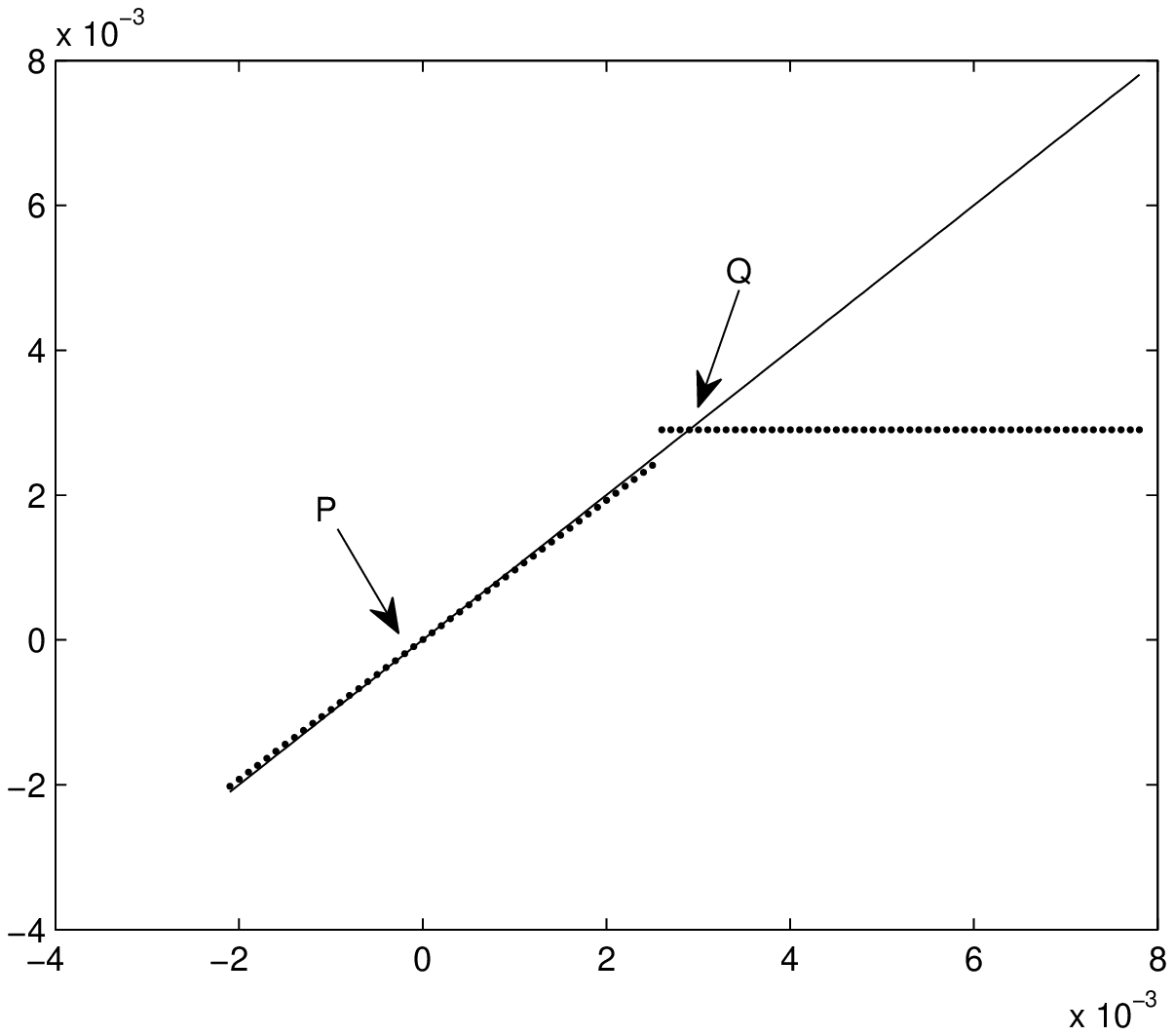}
\begin{tabular}{lr}
\hskip-25pt\includegraphics[width=.56\textwidth]{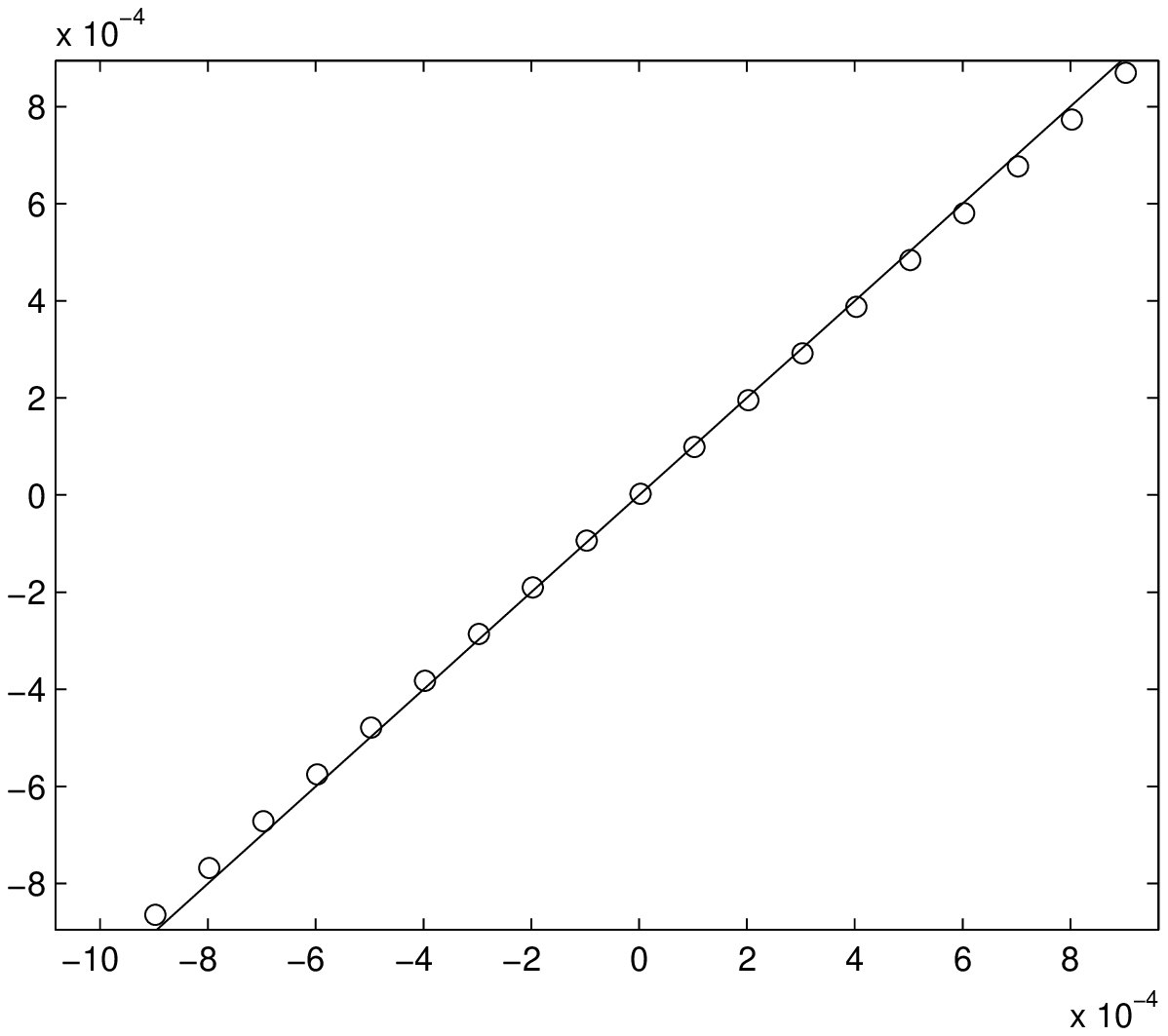} 
\hskip-20pt\includegraphics[width=.56\textwidth]{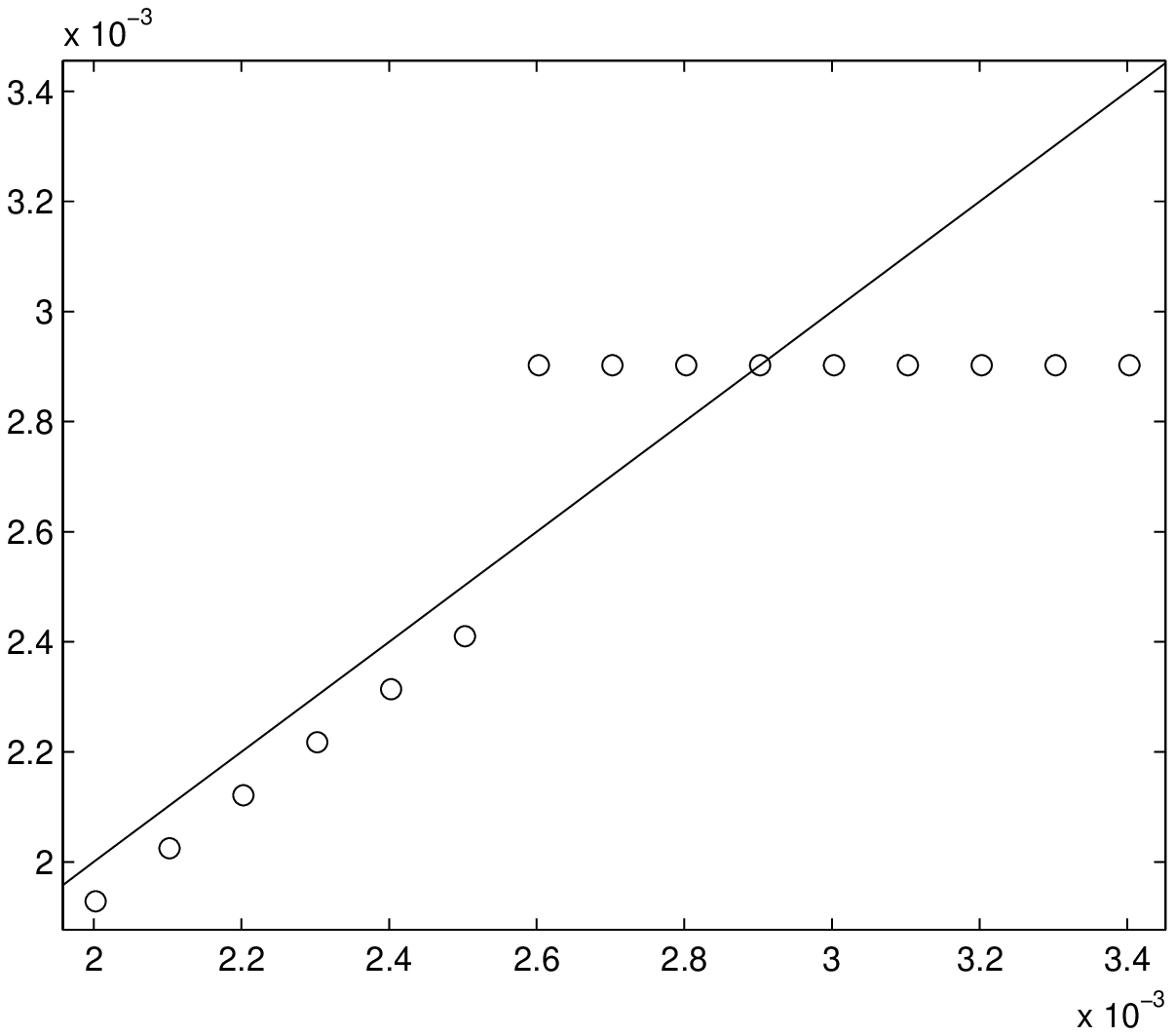}
\end{tabular}
\caption{Test 4. First line: the identity function (solid line) and one of the component of $F$ as a function of one of its argument (dots and line). The function turns to be piecewise contractive with two fixed points $P$ and $Q$. Second line: zoom around $P$ (left) and $Q$ (right).}
\label{T4_contraz_a_tratti}
\end{figure}

\noindent\hskip-3pt{\bf Test 3}\\
Here we present a two-dimensional test with an easy coupled dynamics. We choose   
\begin{equation}\label{dyn2Deasy}
f(x,y,a_1,a_2)=(a_2,a_1)\,,
\end{equation}
with cost functions
\begin{equation}\label{cost2Deasy}
\psi_{i}(x,y,a_1,a_2)=
\left\{
\begin{array}{ll}
\sqrt{x^2+y^2} & \mbox{if } \sqrt{x^2+y^2}>1 \\
0 & \mbox{otherwise}
\end{array}
\right.
,\quad i=1,2.
\end{equation}
%\note{Anche se ognuno gioca per se', questo e' un gioco {\em incidentalmente} cooperativo???}
In this game the two players have the same cost function and want to steer the dynamics in the unit ball centred in $(0,0)$ where the cost is 0. Considering the symmetry of the data, we expect $u_1=u_2$. 
The numerical domain is $\Omega=[-2,2]^2$ and it is discretized by $51\times 51$ nodes.
The sets of admissible controls are $A_1=A_2=[-1,1]$, and they are discretized choosing $A_1^\#=A_2^\#=\{-1,0,1\}$. For $i=1,2$ we set the initial guess $U_i^{(0)}$ equal to a large constant. 
Convergence is reached in a few hundreds of iterations, the results are shown in Fig.\ \ref{T3}.\\
The value functions show the expected behaviour, being equal to 0 in the unit ball centred in $(0,0)$ and growing uniformly in every directions.
We also computed numerically the infinity norm of the Jacobian matrix of $F$, finding
$$
\|J_F(\overline U)\|_\infty=\frac{1}{1+h}
$$
(see (\ref{normJFU})) in three cases for $\overline U$: the final solution, the initial guess and an intermediate value of the fixed point algorithm.
\vskip0.2cm
\noindent{\bf Test 4}\\
The last test is devoted to the investigation of a case where the algorithm does not converge. Indeed, choosing
\begin{equation}\label{dyn2Ddifficult}
f(x,y,a_1,a_2)=(a_1+a_2,a_1-a_2)\,,
\end{equation}
\begin{equation}\label{cost2Ddifficult}
\psi_{i}(x,y,a_1,a_2)=x^2+y^2,\quad i=1,2
\end{equation}
and the other parameters as in Test 3, only some values stabilize, whereas others oscillate. In Fig.\ \ref{T4_U} we show the surfaces of $U_1$, $U_2$ obtained after 1000 iterations.
Let us focus our attention on two nodes with different features. In the first node we observe convergence  whereas in the second we observe an oscillating behavior of the approximating sequence.
%We observed that the value of the node corresponding to the central point $(0,0)$ of $U_1$ converges to a specific value. 
Let $j_0\in\{1,...,2N\}$ be the node corresponding to the central point $(0,0)$ for $U_1$. We freeze all the values of $U=(U_1,U_2)$ except $(U)_{j_0}$, which is replaced by a real variable $s$. 
In Fig.\ \ref{T4_contraz_a_tratti} 
we plot the component $F_{j_0}(U)$ of $F(U)$ as a function of $s$, i.e. $F_{j_0}(s)=F_{j_0}((U)_1,...,(U)_{j_0-1},s,(U)_{j_0+1},...,(U)_{2N})$. Compared with the identity function, it is immediately clear that it is discontinuous and piecewise contractive with two fixed points (see Definition \ref{def:pc}), labelled $P$ and $Q$ (see Fig.\ \ref{T4_contraz_a_tratti}). We observe that by our algorithm the value of the node $j_0$ converges to the fixed point $Q$.

\begin{figure}[h!]
\begin{tabular}{lr}
\hskip-25pt\includegraphics[width=.56\textwidth]{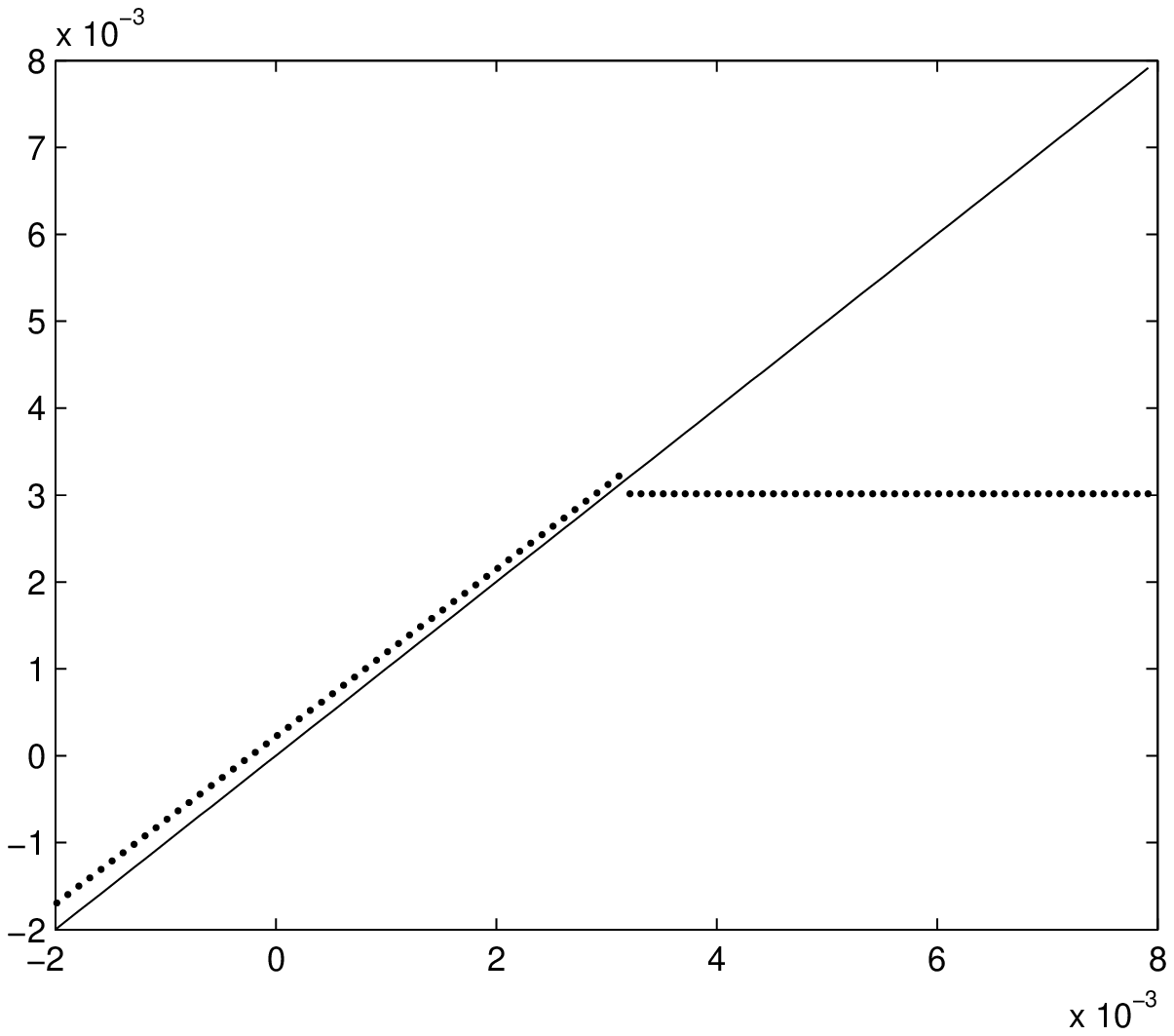}
\hskip-20pt\includegraphics[width=.56\textwidth]{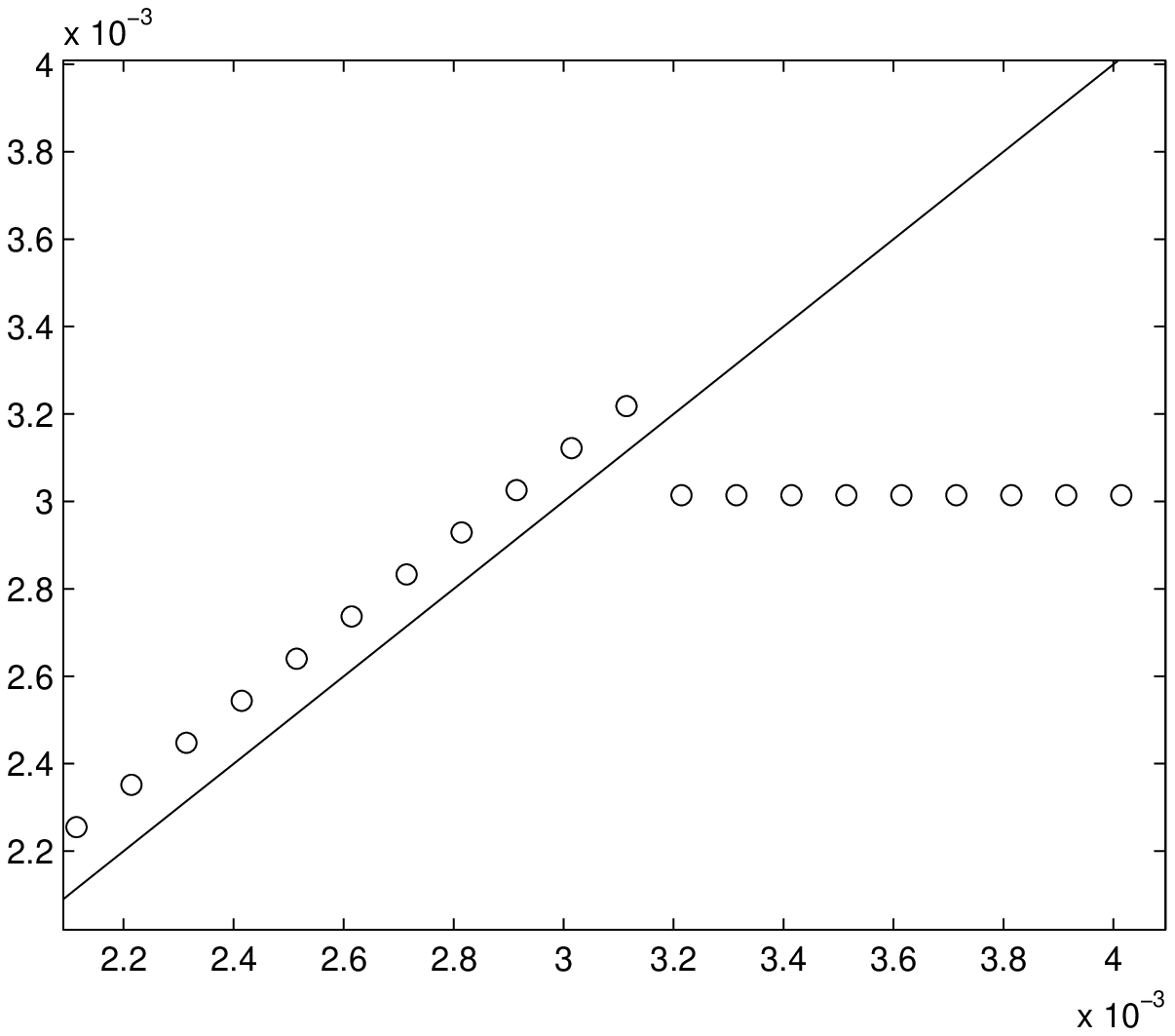} 
\end{tabular}
\caption{Test 4. The identity function (solid line) and one of the component of $F$ as a function of one of its argument (dots and line). The function turns to be piecewise contractive with no fixed points (left). Zoom around the discontinuity (right).}
\label{T4_contraz_a_tratti_bis}
\end{figure}

Conversely, the value of the node corresponding to the point $(0,0+\Delta x)$ does not reach convergence. Similarly as before, we plot the component of the function $F$ corresponding to that point, obtaining the result shown in Fig.\ \ref{T4_contraz_a_tratti_bis}.
In this case the function is piecewise contractive with no fixed points. The value of the node oscillates between two values around the discontinuity, as expected.

\end{document}